**scientific reports**

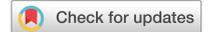

# OPEN    A physical memristor based Muthuswamy–Chua–Ginoux system

Jean-Marc Ginoux[1]✉, Bharathwaj Muthuswamy[2], Riccardo Meucci[3,4], Stefano Euzzor[3], Angelo Di Garbo[5] & Kaliyaperumal Ganesan[6]

In 1976, Leon Chua showed that a thermistor can be modeled as a memristive device. Starting from this statement we designed a circuit that has four circuit elements: a linear passive inductor, a linear passive capacitor, a nonlinear resistor and a thermistor, that is, a nonlinear "locally active" memristor. Thus, the purpose of this work was to use a physical memristor, the thermistor, in a Muthuswamy–Chua chaotic system (circuit) instead of memristor emulators. Such circuit has been modeled by a new three-dimensional autonomous dynamical system exhibiting very particular properties such as the transition from torus breakdown to chaos. Then, mathematical analysis and detailed numerical investigations have enabled to establish that such a transition corresponds to the so-called route to Shilnikov spiral chaos but gives rise to a "double spiral attractor".

Michael Faraday (1791-1867) is generally well known for his contributions to the study of *electromagnetism* and *electrochemistry*. However, according to Orton[1], while investigating the effect of temperature on the conductivity of "sulphuret of silver" (silver sulfide) in 1833 he discovered what is now considered as a *thermistor*[2]. Thermistor, i.e., thermal resistor, is thus an electrical-resistance element made of a semiconducting material the resistance of which varies with temperature. Thermistor production was then difficult, and applications were limited. One century had to pass before thermistors became commonly used by commercial manufacturers. In the 1930s, Samuel Ruben (1900-1988) invented the first commercial thermistor that he called "electrical pyrometer resistance" (Patents N° 2,021,491). He explained that:

> This invention relates to an electrical pyrometer and specifically to one utilizing the resistance change of a metallic compound with heat to indicate temperature changes. The characteristic property of the material employed for the temperature indicating resistance element is one having a high negative resistance coefficient.

In the early 1940s, development of a technique improving the overall consistency and repeatability of many manufacturing processes enabled the large-scale commercialization of thermistor. At that time, thermistor was most frequently used for protection, regulation, and compensating for temperature in electronic circuits. By the 1960s, thermistors were being used in the aerospace industry. Then, John Steinhart and Stanley Hart found a function modeling *thermistor characteristics*, i.e., the resistance according to the temperature which proved to be suitable for a wide variety of thermistors for ranges of a few degrees to a few hundred degrees[3]. In 1976, five years after Chua[4] had postulated a missing circuit element that he called *memristor*, Chua published a paper with Sun Kang[5] in which they recalled that:

> Thermistors have been widely used as a linear resistor whose resistance varies with the ambient temperature. In particular, a *negative-temperature coefficient* thermistor is characterized by

$$v_T = R_0(T_0) \exp\left[\beta\left(\frac{1}{T} - \frac{1}{T_0}\right)\right] i \triangleq R(T)i \qquad (1)$$

> where $\beta$ is the material constant, $T$ is the thermistor temperature and $T_0$ the room temperature both in kelvin. The constant $R_0(T_0)$ denotes the cold temperature resistance at $T = T_0$. The instantaneous tem-perature $T$, however, is known to be a function of the power dissipated in the thermistor and is governed by the heat transfer equation

$$p(T) = v_T(t)i(t) = \delta(T - T_0) + c\frac{dT}{dt} \tag{2}$$

where $c$ is the heat capacitance and $\delta$ is the dissipation constant of the thermistor, defined as the ratio of a change in the power dissipation to the resultant change in the body temperature. Substituting (1) into (2) and by rearranging terms, we obtain

$$\frac{dT}{dt} = -\frac{\delta}{c}(T - T_0) + \frac{R_0(T_0)}{c}\exp\left[\beta\left(\frac{1}{T} - \frac{1}{T_0}\right)\right]i^2 \tag{3}$$

[1]CNRS, CPT, Aix Marseille Univ, Université de Toulon, Marseille, France. [2]Quantum and Nonlinear Engineering Systems, Plainsboro, NJ, USA. [3]Istituto Nazionale di Ottica, Consiglio Nazionale delle Ricerche, Firenze, Italy. [4]Department of Physics and Astronomy, Università di Firenze, Firenze, Italy. [5]Istituto di Biofisica, Consiglio Nazionale delle Ricerche, Pisa, Italy. [6]School of Information Technology and Engineering, Vellore Institute of Technology, Vellore, Tamilnadu 632014, India. ✉email: ginoux@univ-tln.fr

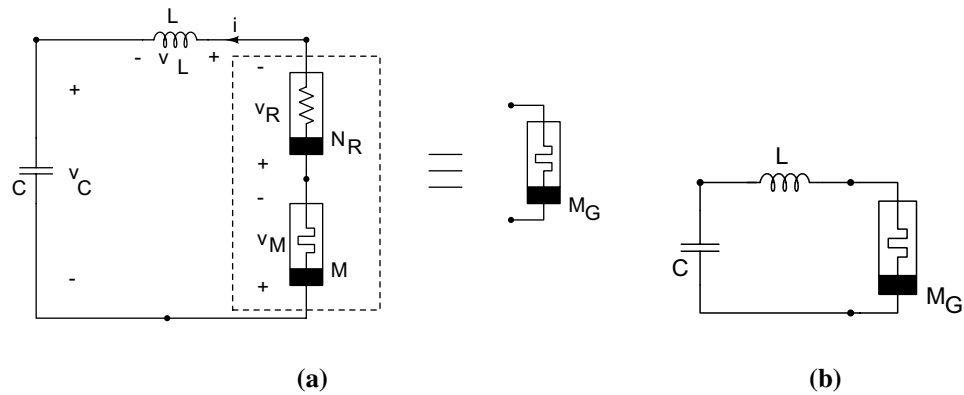

**Figure 1.** (**a**) The Muthuswamy–Chua–Ginoux (MCG) circuit is a generalization of the (**b**) Muthuswamy–Chua circuit, because a "nonlinear resistor in series with a memristor is a memristor" (see Muthuswamy and Banerjee for a proof of this theorem[13]).

We observe from (1) and (3) that *a thermistor* is in fact not a memoryless temperature-dependent linear resistor—as usually assumed—but rather a *first-order time-invariant current-controlled memristive*.

On May 1st, 2008, an ReRAM based memristor was discovered by Strukhov et al.[6]. Their work triggered a renewed interest in memristors and their applications across widely different fields which has consistently grown since. During these last five years, Rajamani et al.[7] analyzed an electronic oscillator circuit designed by connecting an inductor in series with a "locally-active" *Positive Temperature Coefficient* (PTC) memristor and a battery. Then, Sah et al.[8] have considered a "second order memristor which represents the model of a physical device called *Positive Temperature Coefficient* (PTC) and *Negative Temperature Coefficient* (NTC) thermistor connected in series." The aim of this work is to investigate a circuit consisting of a linear passive *capacitor*, a linear passive *inductor*, a nonlinear *resistor* and a *Negative Temperature Coefficient thermistor*, that is, a nonlinear "locally active" volatile *memristor*. Notice that a few implementations of chaotic circuits have recently used ReRAM memristive models[9–11]. However, memristive models for ReRAM devices such as the HP memristor is still the subject of research[12]. In contrast, the memristive model for thermistors is well established and the device can be bought as an off-the-shelf component with a variety of parameter options[13].

### Circuit setup

Starting from the Muthuswamy–Chua circuit[13], we designed a circuit (system) consisting of a linear passive *capacitor* of capacitance $C$, a linear passive *inductor* of inductance $L$, a nonlinear *resistor* $N_R$ and a *thermistor*, that is, a nonlinear "locally active" *memristor* $M$. A schematic is shown in Fig 1.

By Kirchhoff's laws, we have: $v_C + v_L + v_R + v_M = 0$ with $v_C = q/C$, $v_L = L di/dt$, $v_R = f(i)$ and $v_M = R(T)i$. $v_R$ is the voltage across the the nonlinear resistor modeled by a function $f(i)$ defined below. Taking into account the Eq. (3) and since the current intensity across the capacitor is given by $i = dq/dt$, we obtain the following set of equations:

$$\begin{aligned}\frac{dv_C}{dt} &= \frac{i}{C}, \\ \frac{di}{dt} &= -\frac{1}{L}\left[v_C + f(i) + R(T)i\right], \\ \frac{dT}{dt} &= \frac{R(T)}{c}i^2 - \frac{\delta}{c}(T - T_0),\end{aligned} \tag{4}$$

Following Balthasar Van de Pol[14], we model the nonlinear resistor as a cubic function of the current, that is, $f(i) = ai + bi^3$ where $a$ and $b$ are constants. The *thermistor characteristics* (the curve representing the resistance $R$ as a function of the temperature $T$) is modeled using the classical Steinhart–Hart equation ($T^{-1} = A + B \log R + C(\log R)^3$) which is a third-order approximation[3]. By considering that the third Steinhart-Hart coefficient $C = 0$ we obtain the $\beta$-parameter equation (1) which thus corresponds to a second-order approximation. Nevertheless, from both analytical and numerical point of view it is quite difficult to investigate the model (4) with exponential functions. Therefore, physically, we can restrict the thermistor to a small temperature change and hence approximate $R(T)$ by its Taylor series expansion up to the second order:

$$R(T) = R_0(T_0) \exp\left[\beta\left(\frac{1}{T} - \frac{1}{T_0}\right)\right] = R_0\left[1 - \frac{\beta}{T_0^2}(T - T_0) + \frac{\beta(\beta + 2T_0)}{2T_0^4}(T - T_0)^2 + O\big((T - T_0)^3\big)\right] \quad (5)$$

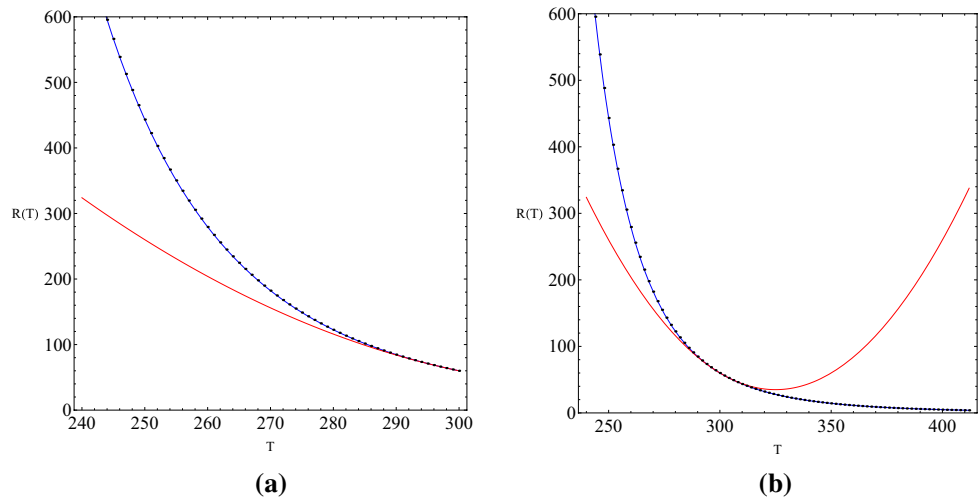

**Figure 2.** *Thermistor characteristics* in blue and its Taylor series expansion in red.

Let us notice that the left hand side of Eq. (5) is a second-order approximation of the Steinhart-Hart equation. As a consequence, its Taylor series expansion must be limited to the second order. Moreover, numerical investigations have shown that if $R(T)$ is approximated by its first order Taylor series expansion, the dynamical system (4) can only exhibit relaxation oscillations and not chaotic behavior. In order to simplify the study of system (4), let's pose for the variables:

$$x = V_C \quad ; \quad y = i \quad ; \quad z = T - T_0.$$

and for the parameters:

$$\alpha = C \quad ; \quad \eta = L \quad ; \quad \theta = \frac{R_0}{c} \quad ; \quad \gamma = -\frac{R_0}{c}\frac{\beta}{T_0^2} \quad ; \quad \mu = \frac{R_0}{c}\frac{\beta(\beta + 2T_0)}{2T_0^4} \quad ; \quad \epsilon = \frac{\delta}{c}.$$

In Fig. 2, we have plotted both *thermistor characteristics* and its Taylor series expansion (5) while using the NTC thermistor B57236S0250M000 (This device can be purchased and its specifications can be found from the following website: https://product.tdk.com/info/en/documents/data_sheet/50/db/icl_16/S236.pdf) the parameters of which are the following: $R_0 = 60\Omega$ and $\beta = 3000$ K.

Then, we have computed the coefficient of determination $\mathcal{R}^2$ for $T \in [240, 300]$ to measure how well the Taylor series expansion fits with the *thermistor characteristics* and found $\mathcal{R}^2 = 0.9726$ ($p$ value $< 0.025$) which indicates a quite good fit[15].

Thus, we obtain what we call the Muthuswamy–Chua–Ginoux (MCG) system:

$$\begin{aligned}\frac{dx}{dt} &= \frac{y}{\alpha}, \\ \frac{dy}{dt} &= -\frac{1}{\eta}[x + f(y) + R(z)y], \\ \frac{dz}{dt} &= R(z)y^2 - \epsilon z,\end{aligned} \quad (6)$$

where $f(y) = ay + by^3$ and $R(z) = \mu z^2 + \gamma z + \theta$.

In order to justify the novelty of our system, we apply the criteria for publication standard for chaotic systems set forth by Sprott[16].

1. The system should credibly model some important unsolved problem in nature and shed insight on that problem.

Ever since the announcement of a memristor by Hewlett-Packard labs[6], memristor based chaotic circuits abound[17]. Nevertheless, such circuits do not use readily available physical memristors such as thermistors, discharge tubes or junction diodes[13]. Our work hence resolves the unsolved question: "can a physical memristor based electronic circuit exhibit chaos"? Moreover, the insight gained is that a nonlinear resistor in series with a memristor (still a memristor, refer to Fig. 1) can be used to appropriately *shape* the nonlinear memristance for chaotic behaviour.

2. The system should exhibit some behavior previously unobserved.

    We observe, for the very first time in a memristor based chaotic circuit, the Shilnikov spiral chaos[18] but with exotic phenomenon such as "multiple" *reverse period-doubling cascade*.

3. The system should be simpler than all other known examples exhibiting the observed behavior.

    Referring to Fig. 1, one can see that we only have four elements in series. Analog memristor based chaotic circuits that use emulators require quite a lot more components[17]. Although this paper uses a cubic characteristic for the nonlinear resistor, one could simply use a piecewise-linear characteristic and hence implement the nonlinear resistor in Fig. 1 using a single op-amp[17]. Please refer to the appendix for electronic circuit implementation.

Let's notice that the flow of this system (6) is invariant under the symmetry $(-x, -y, z) \to (x, y, z)$. Hence if $(x(t), y(t), z(t))$ is a solution of system (6), so is $(-x(t), -y(t), z(t))$.

**Parameters set.** From Eq. (5), it follows that its right hand side is strictly positive since it is an exponential function modeling the resistance of the thermistor which is supposed to be positive. As a consequence, the same should be true for its left hand side, i.e., its Taylor series expansion. Thus, the second degree polynomial $R(z) = \mu z^2 + \gamma z + \theta$ must be positive. This implies that its discriminant $\gamma^2 - 4\mu\theta$ must be negative and that $\mu$ must be positive. Otherwise, the resistance could take negative values which is impossible from a thermodynamics point of view. Moreover, since $\theta = R_0/c$ it must also be positive. Notice that the negativity of $\theta$ would imply the positivity of the discriminant which is here precluded. Such considerations lead to the following conditions:

$$\mu > 0 \quad ; \quad -2\sqrt{\mu\theta} < \gamma < 2\sqrt{\mu\theta} \quad ; \quad \theta > 0.$$

So, first for the stability analysis of the dynamical system (6) presented in the next section, we will use the following parameters set which meet the previous conditions and which facilitate the identification of the transition from torus breakdown to "double spiral chaos":

$$a = -6, b = 3, \eta = 12.2, \mu = 3, \gamma = -2, \theta = 3, \epsilon = 0.6.$$

Then, in the Appendix, dynamical system (6) will be recast in dimensionless form which is more suitable for analog simulations with a parameters set including the NTC thermistor B57236S0250M000 features ($R_0 = 60\,\Omega$, $\beta = 3000$ K). These analog simulations could be thus used to highlight such a transition.

## Stability analysis

**Fixed points.** Fixed points are determined while using the classical nullclines method. MCG system (6) has the origin $O(0, 0, 0)$ as the unique fixed point.

**Jacobian matrix and eigenvalues.** The Jacobian matrix of MCG system (6) reads:

$$J = \begin{pmatrix} 0 & \frac{1}{\alpha} & 0 \\ -\frac{1}{\eta} & -\frac{1}{\eta}(a + 3by^2 + \mu z^2 + \gamma z + \theta) & -\frac{y}{\eta}(2\mu z + \gamma) \\ 0 & 2(\mu z^2 + \gamma z + \theta)y & (2\mu z + \gamma)y^2 - \epsilon \end{pmatrix} \quad (7)$$

By replacing the coordinate of the fixed point, i.e., the origin in the Jacobian matrix (7) one obtains the Cayley-Hamilton third degree eigenpolynomial which can be easily factorized as follows:

$$(\lambda + \epsilon)\left[\alpha\eta\lambda^2 + \alpha(a + \theta)\lambda + 1\right] = 0 \quad (8)$$

Thus, there is one real eigenvalue:

$$\lambda_1 = -\epsilon \quad (9)$$

which is negative since the parameter $\epsilon$ is positive and two eigenvalues:

$$\lambda_{2,3} = \frac{1}{2\eta}\left[-(a + \theta) \pm \sqrt{(a + \theta)^2 - 4\frac{\eta}{\alpha}}\right]. \quad (10)$$

**Eigenvalues' properties.** The value $\alpha(a + \theta)^2 - 4\eta/\alpha$ depends on $\alpha$ which we will be chosen as *bifurcation parameter*. Since we have posed $\alpha = C > 0$, this value is negative provided that: $0 < \alpha \leqslant 4\eta/(a + \theta)^2$. In this stability analysis, we use the following parameters set: $a = -6, \eta = 12.2, \mu = 3, \gamma = -2, \theta = 3, \epsilon = 0.6$ and $b = 3$. So, we have $a + \theta < 0$. Thus, according to Poincaré[19], the real parts of the eigenvalues (10) are positive and the fixed point, i.e., the origin is a *saddle node* or a *saddle focus*. More specifically, if $0 < \alpha \leqslant 4\eta/(a+\theta)^2$, it is a *saddle focus*. These results are summarized in the Table 1.

| $\alpha$ | 0 | $\frac{4\eta}{(a+\theta)^2} = 5.42$ |
|---|---|---|
| $\Delta$ | − | + |
| $\lambda_1$ | − | − |
| $Re(\lambda_2)$ | + | + |
| $Re(\lambda_3)$ | + | + |
| Nature | Saddle-focus | Saddle-node |

**Table 1.** Fixed point stability.

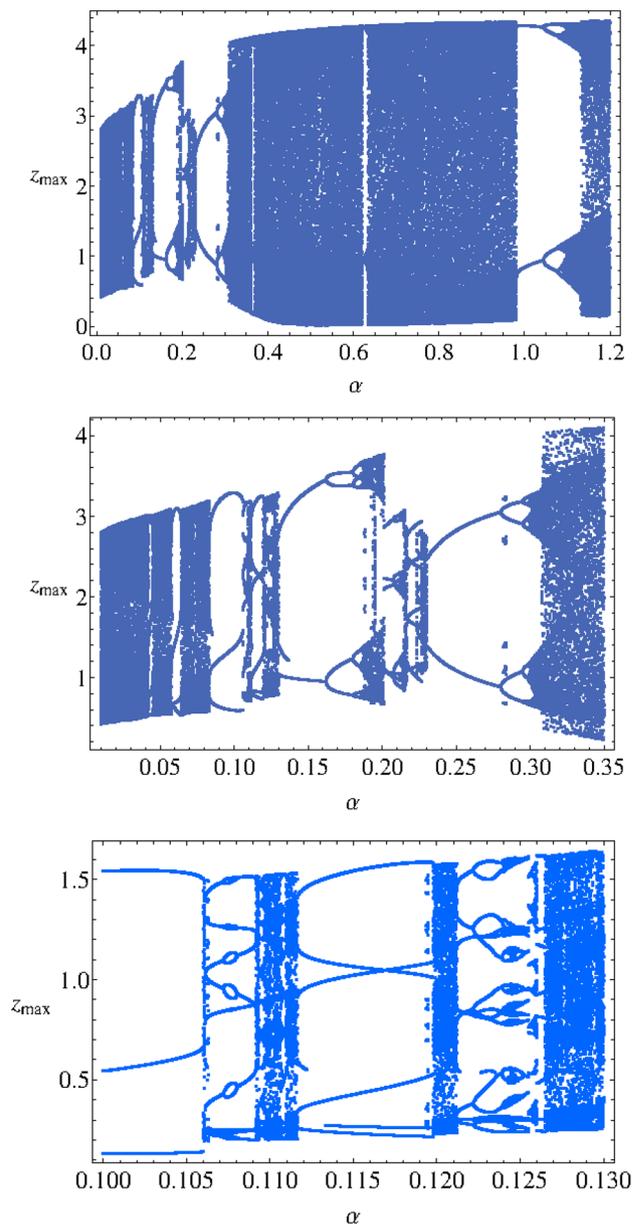

**Figure 3.** Bifurcation diagram $z_{max}$ as function of $\alpha$.

To investigate the effects of the changes of parameter $\alpha$ on the dynamics of the MCG system, a bifurcation diagram has been plotted (see Fig. 3). Then the information can be used to have a better understanding of the phase space orbits plotted in Fig. 4. According to I.M. Ovsyannikov and L.P. Shilnikov[20]:

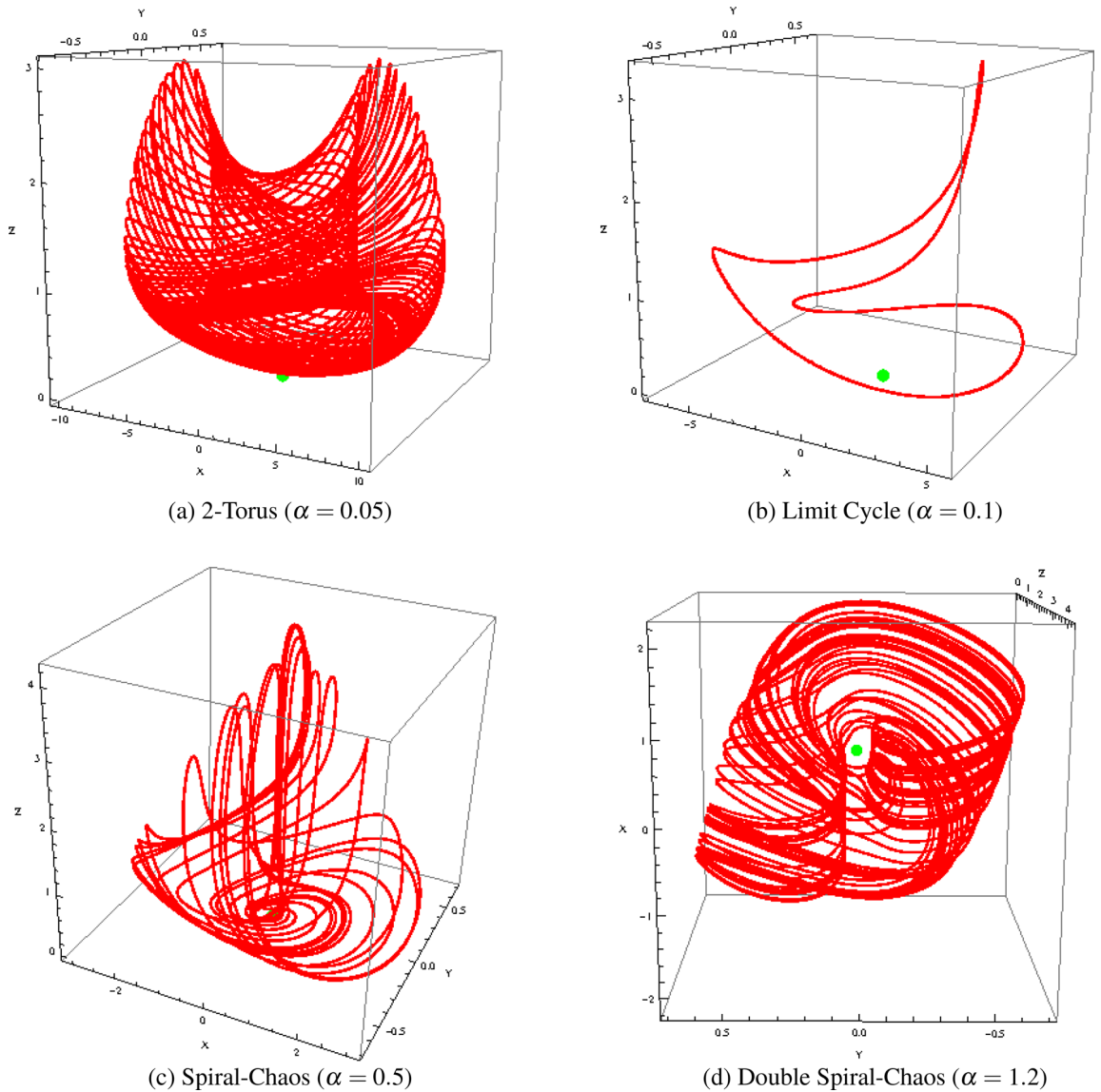

**Figure 4.** Phase portraits of Muthuswamy–Chua–Ginoux system (6) in the phase space for $\alpha = 0.05, 0.1, 0.5$ and 1.2.

The transition to spiral chaos along the lines of the above scenario is usually preceded by either a cascade of period doubling bifurcations (if three-dimensional volumes are contracting) or the collapse of a two-dimensional torus arising from L (a stable periodic motion).

We observe, on the bifurcation diagram (Fig. 3), a "double" *period-doubling cascade* for $\alpha \in [0.001, 0.3]$ and $\alpha \in [0.9, 1.2]$. This "doubling" is due to the fact that, for some particular values of the parameters set, the trajectory curve take the shape of a "double spiral attractor" (see Fig. 4d). Hence, each spiral can perform a *period-doubling cascade*. That's the reason why the bifurcation diagram exhibits two branches for the *period-doubling cascade*.

The transition from torus breakdown to "double spiral chaos" has been highlighted by plotting the phase portrait of MCG system (6) for various values of $\alpha$ (see Fig. 4) and by comparing it with the bifurcation diagram (see Fig. 3). We observe that for $\alpha = 0.05$ (all other parameters are the same as above), the attractor is a *torus* (see Fig. 4a). As $\alpha$ increases between 0.08 and 0.1, a torus breakdown is observed and a *limit cycle* appears (see Fig. 4b). When $\alpha$ increases again between 0.11 and 0.13, a "spiral attractor" appears. Then, for $\alpha \in [0.14, 0.19]$ a *period-doubling cascade* occurs. For $\alpha = 0.2$ a "spiral chaos" attractor is observed while for $\alpha \in [0.21, 0.23]$, a second *period-doubling cascade* occurs. For $\alpha \in [0.24, 0.28]$ a *limit cycle* appears and for $\alpha \in [0.29, 0.30]$, a third *period-doubling cascade* occurs (see Fig. 3). For $\alpha \in [0.31, 0.98]$, a "spiral chaos" attractor is again observed (see Fig. 4c). For $\alpha = 1$ a *limit cycle* appears. Finally, after a short fourth *period-doubling cascade* for $\alpha = 1.1$, a "double spiral attractor" is observed for $\alpha = 1.2$ (see Fig. 4d). The corresponding time series $x(t)$, $y(t)$ and $z(t)$ of Muthuswamy–Chua–Ginoux system (6) have been plotted in Fig. 5 for $\alpha = 0.5$ and 1.2.

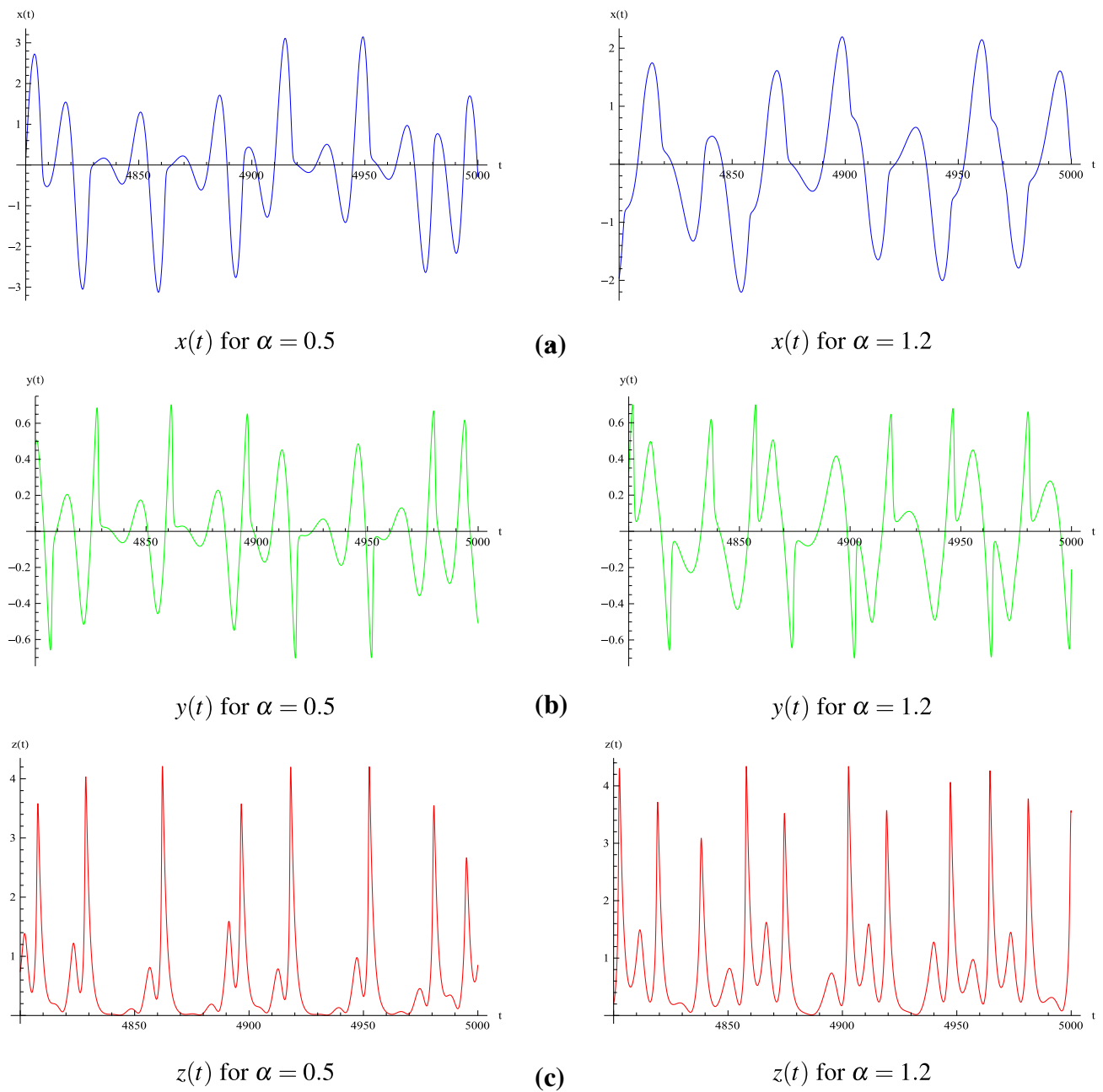

**Figure 5.** Time series of Muthuswamy–Chua–Ginoux system (6) for $\alpha = 0.5$ and 1.2.

In order to confirm such scenario, Lyapunov Characteristic Exponents (LCE) have been computed in each case. The algorithm developed by Sandri[21] for Mathematica has been used to perform the numerical calculation of the Lyapunov characteristics exponents (LCE) of the Muthuswamy–Chua–Ginoux system (6) in each case. LCEs values have been computed within each considered interval ($\alpha \in [0.001, 0.3]$ and $\alpha \in [0.9, 1.2]$). As an example, for $\alpha = 0.05, 0.1, 0.5$ and $1.2$, Sandri's algorithm has provided respectively the following LCEs $(0, 0, -0.13)$, $(0, -0.08, -0.08)$, $(0.08, 0, -0.4)$ and $(0.073, 0, -0.38)$. Then, following the works of Klein and Baier[22], a classification of (autonomous) continuous-time attractors of dynamical system (6) on the basis of their Lyapunov spectrum, together with their Hausdorff dimension is presented in Table 2. LCEs values have been also computed with the Lyapunov Exponents Toolbox (LET) developed by Siu for MatLab and involving the two algorithms proposed by Wolf et al.[23] and Eckmann and Ruelle[24] (see https://fr.mathworks.com/matlabcentral/fileexchange/233-let). Results obtained by both algorithms are consistent.

On Fig. 6, the bifurcation diagram has been plotted for $\alpha \in [0.001, 12]$. Starting from $\alpha = 3$, we observe that a "multiple" *reverse period-doubling cascade* with 1, 2, 3, 4, 5, …branches occurs. Such a number corresponds to the period of each spiral of the "double spiral attractor". As an example, for $\alpha = 7.5$, there are three bifurcation points and so, each spiral has a period 3 while for $\alpha = 11$, there are five bifurcation points and so, each spiral has a period 5 (see Fig. 7).The time series $y(t)$ and $z(t)$ of of Muthuswamy–Chua–Ginoux system (6) plotted in Fig. 8 for $\alpha = 7.5$ and 11 highlights period three and five.

| m | LCE spectrum | Dynamics of the attractor | Hausdorff dimension |
|---|---|---|---|
| $0.001 < \alpha < 0.08$ | $(0, 0, -)$ | 2-Torus | $D = 2$ |
| $0.08 < \alpha < 0.1$ | $(0, -, -)$ | 1−Periodic motion (limit cycle) | $D = 1$ |
| $0.11 < \alpha < 0.13$ | $(+, 0, -)$ | Spiral-chaos | $D = 2.18$ |
| $0.14 < \alpha < 0.19$ | $(0, -, -)$ | $n$-Periodic motion | $D = 1$ |
| $\alpha = 0.20$ | $(+, 0, -)$ | Spiral-chaos | $D = 2.14$ |
| $0.21 < \alpha < 0.23$ | $(0, -, -)$ | $n$-Periodic motion | $D = 1$ |
| $0.24 < \alpha < 0.28$ | $(0, -, -)$ | 1−Periodic motion (limit cycle) | $D = 1$ |
| $0.29 < \alpha < 0.30$ | $(0, -, -)$ | $n$-Periodic motion | $D = 1$ |
| $0.30 < \alpha < 0.98$ | $(+, 0, -)$ | Spiral-chaos | $D = 2.21$ |
| $1 < \alpha < 1.1$ | $(0, -, -)$ | 1−Periodic motion (limit cycle) | $D = 1$ |
| $1.2 < \alpha < 2.9$ | $(+, 0, -)$ | Double spiral-chaos | $D = 2.19$ |

**Table 2.** Lyapunov characteristics exponents of dynamical system (6) for various values of $\alpha$.

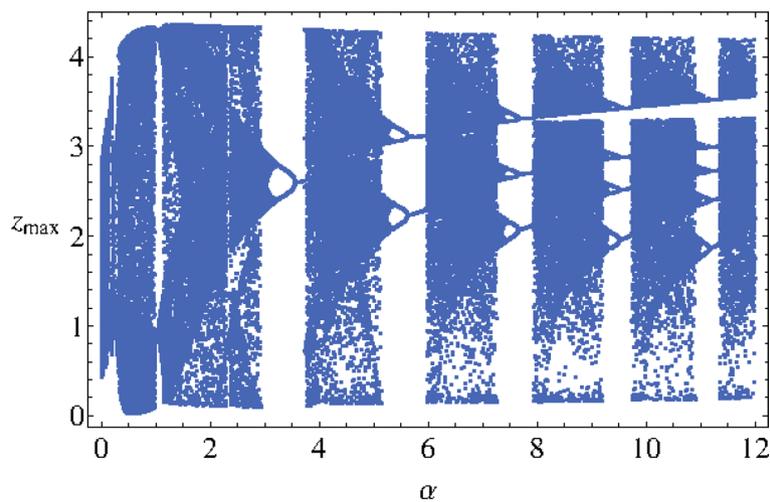

**Figure 6.** Bifurcation diagram $z_{max}$ as function of $\alpha$.

For $\alpha \in [2.9, 12]$, the bifurcation diagram (see Fig. 6) consists in intervals within which the trajectory curve takes either the shape of a "double spiral attractor" or a "$n$-periodic limit cycle".

## Discussion

By using a four element circuit comprising a linear passive inductor, a linear passive capacitor, a nonlinear resistor and a thermistor, that is, a nonlinear "locally active" memristor, we designed a new three-dimensional autonomous dynamical system exhibiting very particular properties such as a transition from a *limit cycle* to *double spiral chaos*. Then, mathematical analysis and detailed numerical investigations have enabled to establish that such a transition corresponds to the route to *Shilnikov spiral chaos* originally described by L.P. Shilnikov[25] and then by I.M. Ovsyannikov and L.P. Shilnikov[26] (see also I.M. Ovsyannikov and L.P. Shilnikov[20] and Shilnikov et al. [27]) but gives rise to a "double spiral attractor". Another scenario leading to spiral chaos, as well as double spiral chaos, and also triple spiral chaos has been proposed by B. Deng[18] twenty-five years ago. Nevertheless, if the attractors he obtained seem to be *topologically equivalent* to those exhibited by the Muthuswamy–Chua–Ginoux system (6) for various values of the *bifurcation parameter $\alpha$*, the route leading to these spiral and double spiral chaotic attractors is different from the "Z-switches" described by Deng[18] but corresponds to the route to *Shilnikov spiral chaos*[25] as highlighted by the *bifurcation diagrams* (see Figs. 3, 6) and the computation of Lyapunov Characteristic Exponents (see Table 2). Double scroll attractor (not double spiral attractor) has been used during the 1990s in cryptography[28]. In 2017, this idea has been patented by Rainer Plaga, Ralph Breithaupt, Sven Freud & Stephan Gieseler, "Chaotic circuit having variable dynamic states as secure information memory," (https://patents.google.com/patent/WO2017097909A1/en). So, it seems possible to imagine that the double spiral attractor could be used for the same application.

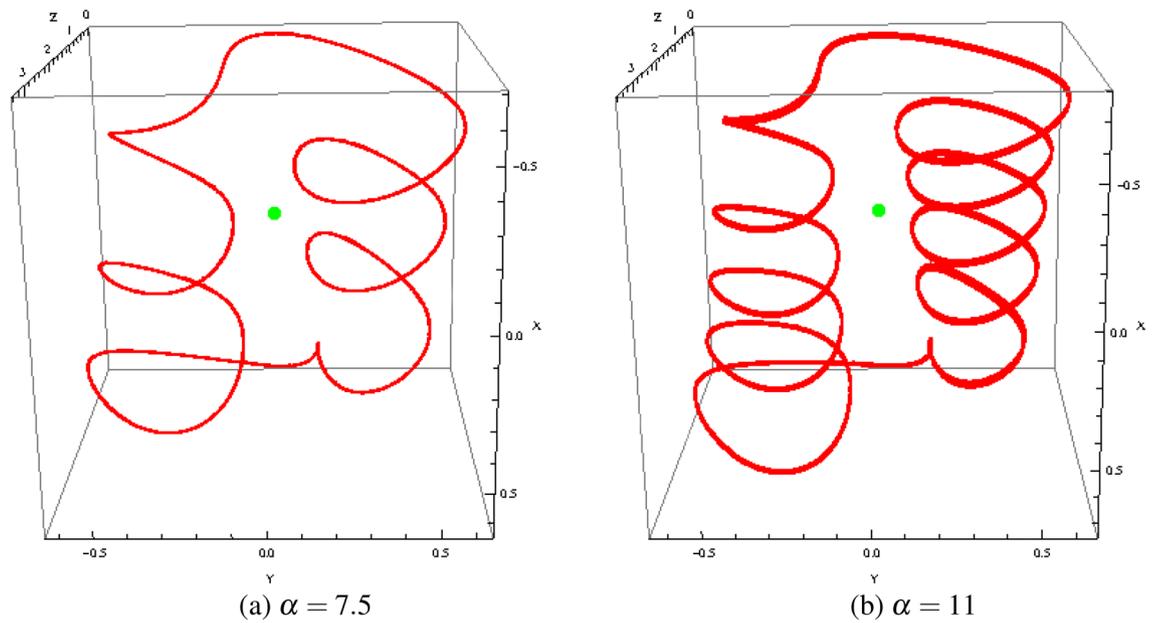

Figure 7. Phase portraits of Muthuswamy–Chua–Ginoux system (6) in the phase space for various values $\alpha$.

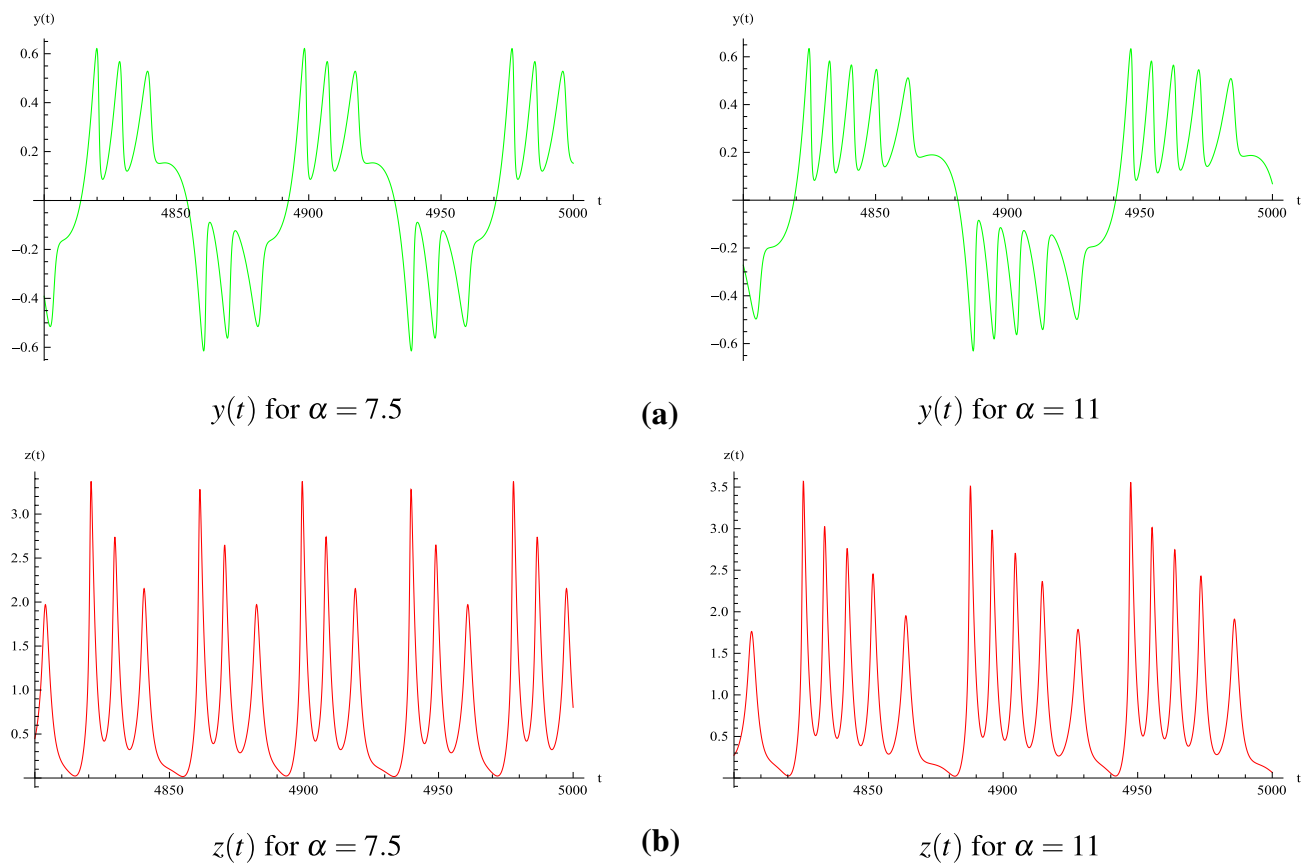

Figure 8. Time series of Muthuswamy–Chua–Ginoux system (6) for $\alpha = 7.5$ and 11.

### Acknowledgements
Authors would to thank Mr. Emmanuel Matte (Distribution Manager, TDK Electronics France SAS), Dr. Roomila Naeck and the reviewers for their helpful advice and comments. Dr. Muthuswamy would like to thank his immediate family (Mr. M. G. Muthuswamy and Mrs. Chandra Muthuswamy) for financial support towards this project.

### Author contributions
All authors have equally contributed to the manuscript.

### Competing interests
The authors declare no competing interests.

### Additional information
**Supplementary information** is available for this paper at https://doi.org/10.1038/s41598-020-76108-z.

**Correspondence** and requests for materials should be addressed to J.-M.G.

**Reprints and permissions information** is available at www.nature.com/reprints.

**Publisher's note** Springer Nature remains neutral with regard to jurisdictional claims in published maps and institutional affiliations.